\documentstyle[psfig,amssymb]{article}
\title{Contact Structures and\\Periodic Fundamental Groups}
\author{Hansj\"org Geiges and Charles B. Thomas}
\date{}

\newcommand{\SO}{\mbox{\rm SO}}

\newcommand{\U}{\mbox{\rm U}}

\newcommand{\Hom}{\mbox{\rm Hom}}
\newcommand{\Ext}{\mbox{\rm Ext}}
\newcommand{\Spin}{\mbox{\scriptsize\rm Spin}}
\newcommand{\Cont}{\mbox{\rm Cont}_5(\pi )}
\newcommand{\cont}{\mbox{\rm Cont}_5}
\newcommand{\om}{\Omega_5^{\Spin}}
\newcommand{\Om}{\Omega_5^{\Spin}(B\pi )}
\begin{document}
\newtheorem{thm}{Theorem}
\newtheorem{prop}[thm]{Proposition}
\newtheorem{lem}[thm]{Lemma}
\newtheorem{cor}[thm]{Corollary}
\newtheorem{defn}[thm]{Definition}
\maketitle
\section{Introduction}
This paper is concerned with the existence of contact structures on
(connected, closed, orientable) 5-manifolds with certain finite fundamental
groups. As such, it constitutes a sequel to~\cite{geig97} (which gave
corresponding existence results for highly connected manifolds of arbitrary
(odd) dimension and some {\em ad hoc} results for finite fundamental groups)
and our joint paper~\cite{geth98}, where we showed that
every 5-manifold $M$ with fundamental group $\pi_1(M)={\Bbb Z}_2$ and
universal cover $\widetilde{M}$ a spin manifold can be obtained from
one of ten `model manifolds' by surgery along a link of 2-spheres and,
as an application of this structure theorem, that every manifold of
this kind admits a contact structure.

In the present paper we combine the ideas of~\cite{geth98} with those
of the extensive literature on the existence of positive scalar
curvature (psc) metrics -- in particular \cite{kwsc90,miya84,rose83,rose86}
(see also \cite{gilk97} and \cite{rost94} for more recent surveys on this
literature) -- to arrive at the following existence result.

\begin{thm}
Let $\pi$ be a finite group of odd order~$|\pi |$ and finite cohomological
period. Furthermore, assume that $|\pi |$ is not divisible by~$9$.
Then every closed $5$-dimensional spin manifold $M$ with fundamental group
$\pi_1(M)\cong \pi$ admits a contact structure.
\end{thm}

Finite groups of odd order and finite cohomological period are
metacyclic with presentation
\[ \left\{ x,y|\, x^m=y^n=1,\, yxy^{-1}=x^r,\, \gcd ((r-1)n,m)=1,
\, r^n\equiv 1\;\mbox{\rm mod}\; m\right\} .\]
This is an old result of Burnside, and it includes the class
of cyclic groups ($m=1$). Geometrically, these groups
are characterized (among groups of odd order) by the property of
acting freely and smoothly on some homotopy sphere.
The case when $n=3$ and $r$ is a primitive cube root of 
1 mod~$m$ is of special interest for 5-manifold topology, because
the corresponding groups
act freely (but not linearly) on~$S^5$. These groups were discussed
in~\cite{geig97}, and Theorem~27 of that paper is a corollary of
Theorem~1 above. More about the connection with that previous paper will
be said at the end of Section~6, where we take a more
geometric view at some of the algebraic arguments in Sections 4 and~5.

The assumption that $|\pi |$ be odd and not divisible by 9 would seem to be
a defect of the proof rather than a defect of nature. In fact, a large
portion of the theory developed in this paper extends to arbitrary groups
of finite cohomological period.  But the case when $|\pi |$ is even 
(which, as regards the existence of psc metrics, has been tackled
successfully in~\cite{bgs97}) seems to present difficulties
of a different order.

For further motivation and some historical comments we refer the reader to
the introduction of~\cite{geig97}. To avoid undue repetition, we assume
the reader to be familiar with the basic definitions of contact
geometry and the fundamental results of contact surgery, due to
Eliashberg~\cite{elia90} and Weinstein~\cite{wein91},
as expounded in Sections 2 and 3 of \cite{geig97}
or the corresponding sections of~\cite{geth98}. On the other hand, while
(equivariant) cobordism arguments have become standard fare in the
literature on psc metrics, this is certainly only true to a much smaller
extent in the contact geometric world, so to make this paper reasonably
self-contained we have chosen to include some arguments which for anyone
familiar with the cited references will certainly cause a sensation
of {\em d\'ej\`a vu}. We shall allow ourselves, however,
to quote liberally from the standard treatise
on periodic maps by Conner and Floyd~\cite{cofl64}.

We now briefly recall the main features of contact surgery on which we
shall rely later on (all details can be found in the beginning sections
of our two earlier papers). In particular, we wish to emphasize the
minor but nonetheless important differences when comparing this with the
surgical arguments for manifolds with psc metrics. Whereas for
the latter any surgery is permitted up to codimension~3, the
restrictions on contact surgeries are:
\begin{itemize}
\item Contact surgery is only possible up to the middle dimension.
\item The sphere along which surgery is performed has to be isotropic, i.e.\
tangent to the contact structure, and it must have trivial conformal
symplectic normal bundle.
\item The framing of the surgery is fixed up to a change in trivializing
the conformal symplectic normal bundle ($CSN$).
\end{itemize}

Because of our restriction to dimension~5, the first point does not entail any
differences between the two theories (contact structures or psc metrics).
The second condition is controlled by an $h$-principle and can be
guaranteed by requiring the given contact structure to have first
Chern class $c_1$ evaluating to zero on 2-spheres.
As regards the third condition, for surgeries along
1-spheres the rank of the $CSN$ is high enough to allow the realization
of any topologically possible framing, and for 2-surgeries we have no
choice of framing because of $\pi_2(\SO_3)=0$. This may serve as an
indication that corresponding existence results for contact structures
in higher dimensions will be much harder to come by.
\section{Periodic fundamental groups}
A finite group $\pi$ is said to have {\em periodic cohomology}
(or simply to be {\em periodic}) if there is some $d>0$ such that
$H^n(\pi )\cong H^{n+d}(\pi )$ for all $n>0$, and the least such $d$
is called the {\em period} of~$\pi$.

We shall use two well-known facts about periodic groups~$\pi$
(cf.~\cite[VI.9]{brow82}):
\begin{itemize}
\item[(1)] Each Sylow subgroup of $\pi$ is cyclic or a generalized quaternion
group (so only the former happens if $|\pi |$ is odd). Indeed, this
statement is equivalent to $\pi$ having periodic cohomology.
\item[(2)] $H_2(\pi )=0$ (in fact, $H_n(\pi )=0$ for $n$ even, $n\geq 2$).
\end{itemize}

Recall from \cite{geig97} that a contact structure $\xi=\ker\alpha$
on a 5-manifold $M$ (where $\alpha$ is a 1-form with $\alpha\wedge (d\alpha
)^2\neq 0$) induces a reduction of the structure group of the
tangent bundle $TM$ to $\U (2)\times 1$. On an orientable 5-manifold~$M$
such a reduction exists
if and only if the third integral Stiefel-Whitney class $W_3(M)=\beta
w_2(M)\in H^3(M;{\Bbb Z})$ vanishes (where $\beta$ denotes the Bockstein
operator of the coefficient sequence ${\Bbb Z}\stackrel{2}{\longrightarrow}
{\Bbb Z}\longrightarrow {\Bbb Z}_2$), or equivalently, if the
second Stiefel-Whitney class $w_2(M)\in H^2(M;{\Bbb Z}_2)$ admits an
integral lift $c_1\in H^2(M;{\Bbb Z})$ (Given $\xi$, such an integral lift
is provided by the first Chern class of the
conformally symplectic bundle $(\xi ,d\alpha )\subset TM$).

The following simple observation shows that we need not be concerned with this
topological obstruction if $\widetilde{M}$ is spin and $\pi_1(M)$ periodic
(and $M$ orientable).

\begin{lem}
Let $\pi$ be a group with periodic cohomology and $M$ a manifold with
$\pi_1(M)\cong \pi$ and universal cover $\widetilde{M}$ a spin manifold.
Then $W_3(M)=0$.
\end{lem}

\noindent
{\em Remark.} This lemma is only included for completeness and future
reference (and to indicate what the optimal statement subsuming
Theorem~1 might be). When we restrict attention to fundamental groups
of odd order~$|\pi |$, then
$H_2(\widetilde{M};{\Bbb Z}_2)\rightarrow H_2(M;{\Bbb Z}_2)$
is surjective (because any 2-cycle in $M$ admits a $|\pi |$-fold covering
by a 2-cycle in~$\widetilde{M}$), and hence $w_2(\widetilde{M})=0$
if and only if $w_2(M)=0$. The arguments in the present paper
require that $M$ be spin, but in~\cite{geth98} Theorem~1 is proved for
$\pi_1(M)={\Bbb Z}_2$ under the weaker assumption that $\widetilde{M}$
be spin.

\vspace{2mm}

\noindent {\em Proof.} By a theorem of Hopf (cf.~\cite{brow82}) there
is an exact sequence
\[ \pi_2(M)\longrightarrow H_2(M)\longrightarrow H_2(\pi ),\]
where the first map is the Hurewicz homomorphism. From statement (2)
above we deduce that the Hurewicz homomorphism is surjective for
periodic fundamental groups. Furthermore, the universal covering map
$\widetilde{M}\rightarrow M$ induces an isomorphism $\pi_2(\widetilde{M})
\rightarrow \pi_2(M)$.

Combining this with the assumption $w_2(\widetilde{M})=0$, we find that
$w_2(M)$ maps to zero under the natural homomorphism
\[ H^2(M;{\Bbb Z}_2)\longrightarrow \Hom (H_2(M),{\Bbb Z}_2).\]
Now consider the commutative diagram built from the Bockstein exact sequence
and the universal coefficient theorem:
\[ \begin{array}{ccccc}
\Ext (H_1(M),{\Bbb Z}) & \rightarrowtail & H^2(M;{\Bbb Z}) & & \\
\mbox{\Large $\downarrow$} & & \mbox{\Large $\downarrow$} & & \\
\Ext (H_1(M),{\Bbb Z}_2) & \rightarrowtail & H^2(M;{\Bbb Z}_2) &
        \twoheadrightarrow & \Hom (H_2(M),{\Bbb Z}_2)\\
 & & \,\;\;\mbox{\Large $\downarrow$}\beta & & \\
 & & H^3(M;{\Bbb Z}) & &
\end{array} \]
By the right exactness of $\Ext (G,-)$, the homomorphism between the
$\Ext$ groups in this diagram is surjective. Then a simple diagram
chase allows to conclude that $W_3(M)=\beta w_2(M)=0$. \hfill $\Box$
\section{Contact groups and a reduction theorem}
For any finite group $\pi$ let $\Om$ be the 5-dimensional spin bordism group
of~$\pi$. In other words, elements of this group are equivalence
classes of pairs $(f :V\rightarrow B\pi ,\sigma )$, where $(V,\sigma )$
is a closed 5-dimensional spin manifold with spin structure $\sigma$
and $f$ is a continuous map into the classifying space of~$\pi$,
and spin bordant pairs are
regarded as equivalent. Define $\Cont\subset\Om$ as the set of all classes
with representatives of the form $(f :V\rightarrow B\pi ,\sigma )$,
where $V$ admits a contact structure
defining the orientation given by~$\sigma$ and with first Chern class
$c_1=0$ on the image of $\pi_2(V)$ in~$H_2(V)$.

Changing from a contact structure $\xi =\ker\alpha$ to $\xi =\ker (-\alpha )$,
which amounts to changing the coorientation of $\xi$, changes the orientation
determined by the volume form $\alpha \wedge (d\alpha )^2$. Thus, if $V$
admits a spin and a contact structure, it does so for either
orientation, which allows to take inverses in $\Cont$. The sum operation
in $\Om$ is given by disjoint union, and $\Cont$ always contains the zero
element of $\Om$, represented by $S^5$ and the constant map into
$B\pi$, say, so $\Cont$ is actually a subgroup of~$\Om$.

\begin{thm}
\label{thm:cobord}
Let $(M,\sigma)$ be a connected, closed $5$-dimensional
spin manifold with
fundamental group $\pi$ and let $f:M\rightarrow B\pi$ be the classifying map
of the universal cover $\widetilde{M}\rightarrow M$.
If $(f:M\rightarrow B\pi ,\sigma)$ represents an element in $\Cont$,
then $M$ admits a contact structure.
\end{thm}
The following statement is an immediate consequence of this theorem and the
fact that $\Cont$ always contains the zero element.
\begin{cor}
If $(f:M\rightarrow B\pi ,\sigma )$ as in the theorem represents the
zero element in $\Om$, that is, if $M=\partial W$ with $W$ a compact spin
manifold and $f$ extends over~$W$, then $M$ admits a contact structure.
\end{cor}
 
Because of $\om =0$, this corollary includes the result that every simply
connected 5-dimensional spin manifold admits a contact structure
(see~\cite{geig97} for a stronger theorem in this simply connected case).

In view of Theorem~\ref{thm:cobord} we call $\pi$ a {\em contact group} if
\[ \Cont =\Om. \]
Thus, for contact groups the conclusion of Theorem~1 holds. Conversely,
the result of \cite{geth98} implies that ${\Bbb Z}_2$ is a contact
group, if one observes that any class in $\Omega_5^{\scriptsize\rm Spin}
(B{\Bbb Z}_2)$ can be represented by a manifold with fundamental
group~${\Bbb Z}_2$ (see Section~4 for the corresponding statement
for~${\Bbb Z}_p$, $p$ an odd prime).

It might seem more attractive to require, in the definition of $\Cont$,
that $f$ be the classifying map for the universal cover of~$V$.
Part of the argument for proving Theorem~\ref{thm:cobord} as it stands
could then be used to prove that $\Cont$ is still a subgroup, and the
proof of Theorem~\ref{thm:cobord} with the alternative definition of
$\Cont$ would simplify correspondingly. In some sense, this would be the
approach analogous to the one taken by Rosenberg in~\cite{rose83}.
The present approach is analogous to that of Kwasik and Schultz~\cite{kwsc90}
and has the advantage that we get similar naturality properties for $\Cont$
as they get for a corresponding subgroup $\mbox{\rm Pos}_5(\pi )\subset
\Om$.

Before proving Theorem~\ref{thm:cobord}, we continue with the general set-up
for the proof of Theorem~1.

Given a group homomorphism $h:\pi \rightarrow \pi'$ we have an induced
homomorphism
\[ \begin{array}{rrcl}
(Bh)_*: & \Om & \longrightarrow & \Omega_5^{\scriptsize\rm Spin}(B\pi ')\\
  & \left(f:V\rightarrow B\pi ,\sigma \right) & \longmapsto &
       \left( (Bh)\circ f:V\rightarrow B\pi ',\sigma\right).
\end{array} \]
If $h$ is an inclusion, there is a transfer homomorphism
\[ (Bh)^!:\,\Omega_5^{\scriptsize\rm Spin}(B\pi ') \longrightarrow \Om ,\]
which is defined geometrically as follows: Given
\[ \left( f':V'\rightarrow B\pi',\sigma'\right) \in
\Omega_5^{\scriptsize\rm Spin}(B\pi '),\]
let $\hat{V}\rightarrow V'$ be the principal $\pi'$-bundle defined by~$f'$.
Then the subgroup $h(\pi )\equiv \pi$ of $\pi'$ also acts on $\hat{V}$.
Set $V=\hat{V}/\pi$, let $f:V\rightarrow B\pi$ be the classifying map
of the covering $\hat{V}\rightarrow V$, and lift the spin structure
$\sigma'$ on $V'$ to a spin structure $\sigma$ on $V$ via the covering
$V\rightarrow V'$. Then define
\[ (Bh)^!(f':V'\rightarrow B\pi ',\sigma ')=(f:V\rightarrow B\pi ,\sigma ).\]
We have the following naturality properties of $\Cont$ with respect to these
homomorphisms.

\begin{lem}
\label{lem:naturality}
{\rm (i)} $(Bh)_*$ sends $\Cont$ to $\mbox{\rm Cont}_5(\pi')$.

{\rm (ii)} If $h$ is an inclusion, $(Bh)^!$ sends $\mbox{\rm Cont}_5(\pi')$
to $\Cont$.
\end{lem}

\noindent {\em Proof.}
The first statement is obvious from the construction, and for the second
statement we only need to observe that a contact structure on $V'$ with
$c_1=0$ on 2-spheres lifts to such a structure on~$V$.\hfill $\Box$

\vspace{2mm}

The following reduction theorem is the direct analogue of Proposition~1.5
in~\cite{kwsc90}.

\begin{thm}
\label{thm:reduction}
Let $\pi$ be a finite group of odd order, let $p$ be a prime dividing~$|\pi |$,
and let $j_p:\pi_p\rightarrow \pi$ be the inclusion of a Sylow
$p$-subgroup. Then a class $\alpha\in\Om$ lies in $\Cont$ if and only
if the images $(Bj_p)^!\alpha\in\Omega_5^{\scriptsize\rm Spin}(B\pi_p)$
under the transfer homomorphism of $j_p$ lie in $\mbox{\rm Cont}_5(\pi_p)$
for all~$p$.
\end{thm}

The proof of this theorem can in principle be taken word for word from
the cited paper. For the reader's convenience we reproduce this
proof in Section~6, including additional details of the `standard'
transfer arguments used by Kwasik and Schultz. For our computations in the
subsequent sections we have to discuss the Atiyah-Hirzebruch bordism
spectral sequence, and with details about this spectral sequence at hand the
mentioned transfer arguments become quite transparent.

Using property (1) of periodic groups, we see that it suffices now to
prove Theorem~1 for cyclic groups $\pi={\Bbb Z}_{p^k}$ with
$p$ an odd prime (and $k=1$ only for $p=3$).
With Theorem~3 in mind we see that we are left with showing that
these cyclic groups are contact groups. This will be done in the following
two sections.

\vspace{2mm}

\noindent {\em Proof of Theorem~$3$.}
Write $M_0=M$, $f_0=f$. By assumption, there is a closed (but not
necessarily connected) 5-dimensional
spin manifold $M_1$ admitting a contact structure with $c_1=0$ on
2-spheres, and
a map $f_1:M_1\rightarrow B\pi$ spin bordant to~$f_0$. That is, we have
a 6-dimensional compact spin manifold $W$ (which we may assume
to be connected) with boundary $\partial W=M_1-M_0$,
inducing the given spin structures on $M_0,M_1$, and a map $F:W\rightarrow
B\pi$ restricting to $f_i:M_i\rightarrow B\pi$ on the boundary components.
Write $j_i$ for the inclusion of $M_i$ in $W$ and denote by subscript
`$\#$' induced homomorphisms on homotopy groups. We have the sequence
of homomorphisms
\[ \pi_1(M_0)\stackrel{j_{0\#}}{\longrightarrow} \pi_1(W)
\stackrel{F_{\#}}{\longrightarrow} \pi, \]
where the composition
\[ F_{\#}\circ j_{0\#} = (F\circ j_0)_{\#}=f_{0\#} \]
is an isomorphism by our hypotheses. We thus obtain a split exact sequence
\[ 1 \longrightarrow \ker F_{\#} \longrightarrow \pi_1(W)
\stackrel{F_{\#}}{\longrightarrow} \pi\longrightarrow 1.\]
The group $\ker F_{\#}$ is generated by embedded copies of $S^1$ in
$W$ not meeting the boundary, and performing surgery along these circles will
kill~$\ker F_{\#}$. The choice of framing lies in
$\pi_1(\SO_5)\cong {\Bbb Z}_2$,
and for one of the two framings the surgery will preserve the spin
structure.

So we may assume that $F_{\#}$ and $j_{0\#}$ are isomorphisms. Then the
homotopy exact sequence of the pair $(W,M_0)$ becomes
\[ \pi_2(M_0)\stackrel{j_{0\#}}{\longrightarrow} \pi_2(W)
\longrightarrow \pi_2(W,M_0)\longrightarrow 0.\]
Represent a set of elements of $\pi_2(W)$ generating
\[ \pi_2(W,M_0)\cong \pi_2(W)/j_{0\#}\pi_2(M_0) \]
by smoothly embedded 2-spheres which do
not meet the boundary (which is possible by the Whitney embedding theorem).
Since $W$ is a spin manifold, these spheres have trivial normal
bundle, and surgery along these 2-spheres will kill $\pi_2(W,M_0)$
and preserve fundamental group and spin structure.

We have thus reduced the problem to the case where $(W,M_0)$ is
2-connected. A result of Wall~\cite[Theorem~3]{wall71} says that
homotopical connectivity implies geo\-metrical connectivity in codimension
$\geq 4$, so $(W,M_0)$ is actually {\em geometrically} 2-connected.
This means that $W$, viewed as a cobordism on~$M_0$, contains only handles
of index~$\geq 3$, and thus $M_0$ is obtained from $M_1$ by surgery in
dimension less than or equal to~2.

It remains to be checked that all these surgeries can be performed as
contact surgeries. Clearly there is no problem with 0-surgeries.
The choice of framing of contact 1-surgeries lies in $\pi_1(\U_1)\cong
{\Bbb Z}$ (the conformal symplectic normal bundle of an $S^1$ in
a contact 5-manifold has rank~2). The homomorphism $\pi_1(\U_1)
\rightarrow \pi_1(\SO_4)={\Bbb Z}_2$ induced by inclusion is
surjective, so any topological framing can be realized by
a contact surgery. Furthermore, the framing in $\pi_1(\U_1)$ determines
$c_1$ of the resulting contact manifold, and since all surgeries
preserve the spin structure, we can actually ensure that the
property $c_1|\pi_2=0$ is preserved. Then the remaining surgeries along a link
of 2-spheres can be performed as contact surgeries as well. \hfill $\Box$
\section{Cyclic groups of prime order}
In this section we consider the case $\pi \cong {\Bbb Z}_p$
with $p$ an odd prime.
The fact that all these groups are contact groups is a consequence of
the following proposition, since every lens space $L^5_p$ (indeed, any
quotient of $S^{2n+1}$ under a discrete group acting freely and
linearly, cf.~\cite{geig97}) admits a contact structure, and any such
structure trivially has $c_1|\pi_2=0$, since $\pi_2(L^5_p)=0$.

\begin{prop}
\label{prop:cyclic}
We have $\om (B{\Bbb Z}_3)\cong {\Bbb Z}_9$ and $\om (B{\Bbb Z}_p)
\cong {\Bbb Z}_p\oplus {\Bbb Z}_p$ for $p$ a prime~$\geq 5$. All these groups
are generated by $5$-dimensional lens spaces.
\end{prop}

This proposition is essentially due to Conner and Floyd as far as the
computation of cobordism groups is concerned, and the observation about lens
spaces as generators was made by Rosenberg~\cite{rose86}. We
have not been able to infer this observation from the reference he
quotes, though,  and therefore provide our own proof, which actually
yields a slightly stronger result (see the statement before
Lemma~\ref{lem:numbers}).

\vspace{2mm}

\noindent {\em Proof.}
Write $\Omega_k'$ for $\Omega_5^{\scriptsize\rm Spin}$ or $\Omega_k$ and
$\widetilde{\Omega}_k'(B\pi )$ for the kernel of the homomorphism
\[ \Omega_k'(B\pi )\longrightarrow \Omega_k'(\{ *\} )=\Omega_k' \]
induced by the constant map. Since $\om =0$ we have
$\widetilde{\Omega}_5^{\scriptsize\rm Spin}(B\pi )=\Om$ of course, but
for determining this group it is more convenient to work with reduced bordism
groups.

There is an Atiyah-Hirzebruch spectral sequence for
$\widetilde{\Omega}_*'$ of the form
\[ E^2_{r,s}=\widetilde{H}_r(B\pi ;\Omega_s')\Longrightarrow
\widetilde{\Omega}_*'(B\pi ) \]
(cf.~\cite[Section~7]{cofl64}). We have $\widetilde{H}_r(B{\Bbb Z}_p)\cong
{\Bbb Z}_p$ in positive odd dimensions~$r$ and 0 otherwise, and hence
$\widetilde{H}_*(B{\Bbb Z}_p;{\Bbb Z}_2)=0$. Now $\Omega_*'$ has
only 2-torsion (cf.~\cite{ston68}). In fact, the relevant groups for us
are
\[ \Omega_0'\cong\Omega_4'\cong {\Bbb Z},\;\; \Omega_1^{\scriptsize\rm Spin}
\cong\Omega_2^{\scriptsize\rm Spin}\cong {\Bbb Z}_2,\;\;
\mbox{\rm and}\;\; \Omega_1=\Omega_2=\Omega_3^{\scriptsize\rm Spin}
=0.\]
Thus the spectral sequence collapses and $E^{\infty}_{r,s}=E^2_{r,s}$.
So we obtain the short exact sequence
\[ 0\longrightarrow H_1(B{\Bbb Z}_p)\longrightarrow
\widetilde{\Omega}_5'(B{\Bbb Z}_p)
\stackrel{\mu}{\longrightarrow} H_5(B{\Bbb Z}_p)\longrightarrow 0,\]
that is,
\[ 0\longrightarrow {\Bbb Z}_p\longrightarrow 
\widetilde{\Omega}_5'(B{\Bbb Z}_p)
\longrightarrow {\Bbb Z}_p\longrightarrow 0,\]
where by \cite[(7.2)]{cofl64} the homomorphism $\mu$ is given by
\[ (f:M\rightarrow B{\Bbb Z}_p) \longmapsto f_*[M].\]
The map $\Omega_*^{\scriptsize\rm Spin}\otimes {\Bbb Z}_p\rightarrow
\Omega_*\otimes{\Bbb Z}_p$ given by forgetting the spin structure is
an isomorphism (we only need this in dimension~4, where it follows from
explicit calculations, cf.~\cite{kirb89}). Then the 5-lemma applied to
the two short exact sequences above (for $\Omega_5$ and
$\Omega_5^{\scriptsize\rm Spin}$) shows that
$\widetilde{\Omega}_5^{\scriptsize\rm Spin}(B{\Bbb Z}_p)\rightarrow
\widetilde{\Omega}_5(B{\Bbb Z}_p)$ is an isomorphism (indeed, this
is again true in all dimensions, cf.~\cite{rose83}).

We notice in particular that $\widetilde{\Omega}_5(B{\Bbb Z}_p)$ has
order~$p^2$.
For $p=3$, 5-dimensional lens spaces have order 9 in
$\widetilde{\Omega}_5(B{\Bbb Z}_3)$
according to~\cite[(36.1)]{cofl64}, hence
 $\widetilde{\Omega}_5(B{\Bbb Z}_3)
\cong {\Bbb Z}_9$. For $p\geq 5$, that same theorem states that lens spaces
have order~$p$. So we can define a splitting for $\mu$ by sending
a suitable generator of
$H_5(B{\Bbb Z}_p)$ to the class of some 5-dimensional lens space in
$\widetilde{\Omega}_5(B{\Bbb Z}_p)$, and we see
that $\widetilde{\Omega}_5(B{\Bbb Z}_p)\cong {\Bbb Z}_p\oplus {\Bbb Z}_p$.

In order to prove that $\widetilde{\Omega}_5(B{\Bbb Z}_p)$
is generated by lens spaces
also for $p\geq 5$ we appeal to (34.5) of~\cite{cofl64}, which states
that an element in $\widetilde{\Omega}_5(B{\Bbb Z}_p)$ is zero if and
only if all its mod~$p$ Pontrjagin numbers are zero and thus
implies that it suffices to find two lens spaces $L_p^5$ (for each~$p$)
whose pairs of mod~$p$ Pontrjagin numbers are linearly independent
over~${\Bbb Z}_p$.

We briefly recall the definition of mod~$p$ Pontrjagin numbers,
cf.~\cite[(34.4)]{cofl64}. Choose a generator $d_1$ of
$H^1(B{\Bbb Z}_p;{\Bbb Z}_p)$
and let $d_2\in H^2(B{\Bbb Z}_p;{\Bbb Z}_p)$ be the image of $d_1$ under the
Bockstein operator of the coefficient sequence ${\Bbb Z}
\stackrel{p}{\longrightarrow} {\Bbb Z}\longrightarrow {\Bbb Z}_p$, followed by
mod~$p$ reduction. Then $d_1d_2^2$ is a generator of
$H^5(B{\Bbb Z}_p;{\Bbb Z}_p)$.
Given a 5-dimensional lens space $L_p^5$, let $f:L_p^5\rightarrow
B{\Bbb Z}_p$ be a classifying map for its universal covering.
Specifying a generator for $\pi_1(L_p^5)\cong {\Bbb Z}_p$ amounts to
choosing a homotopy class of classifying maps $f:L_p^5\rightarrow
B{\Bbb Z}_p$. Continue to write $d_i$ for $f^*d_i$, $i=1,2$. Further,
let $p_1\in H^4(L_p^5;{\Bbb Z}_p)$ be the mod~$p$ reduction of the first
Pontrjagin class of~$L^5_p$. Then the mod~$p$ Pontrjagin numbers of
$L^5_p$ are the integers mod~$p$
\[ \beta_0=\langle d_1d_2^2,[L_p^5]\rangle \;\;\mbox{\rm and}
\;\;\beta_1=\langle p_1d_1,[L_p^5]\rangle ,\]
where $[L_p^5]$ is the fundamental cycle of $L^5_p$ and $\langle -,-\rangle$
the Kronecker product. Here $\beta_0$ is always nonzero.

The Pontrjagin classes of lens spaces have been computed by
Folkman~\cite{folk71}, cf.~\cite{miln66}.
For the quotient of $S^5\subset {\Bbb C}^3$ under the
action of ${\Bbb Z}_p$ generated by
\[ T:\, (z_1,z_2,z_3)\longmapsto (\alpha_1z_1,\alpha_2z_2,\alpha_3z_3)\]
with $\alpha_j=\exp (2\pi iq_j/p)$ we have
\[ p_1=(q_1^2+q_2^2+q_3^2)d_2^2 \]
(the choice of a generator $T$ determines $d_1$ and hence~$d_2$).
Replacing $T$ by $T^m$ with $m$ coprime to $p$ amounts to replacing
$q_j$ by $mq_j$ ($j=1,2,3$) and $d_i$ by $kd_i$ ($i=1,2$) with $mk\equiv 1$
mod~$p$. So the mod~$p$ Pontrjagin numbers of $L^5_p(q_1,q_2,q_3)$
modulo the choice of classifying map $f:L_p^5\rightarrow B{\Bbb Z}_p$ are
\[ (\beta_0,\beta_1=(q_1^2+q_2^2+q_3^2)\beta_0)\]
modulo the equivalence relation
\[ (\beta_0,\beta_1)\sim (k^3\beta_0,k\beta_1)\]
for $k$ not divisible by $p$.

The proof of Proposition~\ref{prop:cyclic} is therefore completed with the
following lemma, which proves more than we really need, namely, that it
is possible to find two lens spaces $L^5_{p,1}$ and $L^5_{p,2}$ such
that $[L^5_{p,1},f_1]$ and $[L^5_{p,2},f_2]$ generate $\widetilde{\Omega}_5
(B{\Bbb Z}_p)$ for {\em any} choice of classifying maps $f_i:L^5_{p,i}
\rightarrow B{\Bbb Z}_p$ of their universal coverings. In this lemma we write
\[ Q=q_1^2+q_2^2+q_3^2\;\;\mbox{\rm and}\;\; R=r_1^2+r_2^2+r_3^2.\]

\begin{lem}
\label{lem:numbers}
For any prime $p\geq 5$ there are
triples $(q_1,q_2,q_3)$ and $(r_1,r_2,r_3)$ of
integers mod~$p$ such that the equation
\[ a(k^3\beta_0, kQ\beta_0)+ b(l^3\beta_0', lR\beta_0')
\equiv (0,0)\; \mbox{\rm mod}~p\]
has {\em no} solution $\beta_0,\beta_0',a,b,k,l$ (coprime to~$p$).
\end{lem}

\noindent {\em Proof.}
The pair of equations in the lemma yields
\[ b(Rk^2-Ql^2)l\beta_0'\equiv 0.\]
Since we are assuming $l$ to be coprime to~$p$ we can divide mod~$p$ by
$l^2$ and obtain, by neglecting the factors coprime to $p$ and replacing
$k^2/l^2$ by $k^2$,
\[ Rk^2-Q\equiv 0.\]
We begin with $(q_1,q_2,q_3)=(1,1,1)$ and $(r_1,r_2,r_3)=(1,1,2)$, that is,
$Q=3$ and $R=6$. This yields the equation $6k^2-3\equiv 0$, and hence
$2k^2-1\equiv 0$, since $p\geq 5$. Rewriting this as $2k^2-1=(2n+1)p$
we get
\[ k^2=np+\frac{p+1}{2}.\]
So if $(p+1)/2$ is not a quadratic residue mod~$p$ (e.g.\ if $p=5$),
we are done. Assume, on the contrary, that it is. Then we may take $Q=3$ and 
\[ R\equiv\frac{p+1}{2} +\frac{p+1}{2} + 1^2\equiv 2.\]
This gives $2k^2-3\equiv 0$, and hence
\[ k^2\equiv \frac{p+3}{2} .\]
Again, if $(p+3)/2$ is not a quadratic residue mod~$p$, we are done.
But, since we are assuming that $(p+1)/2$ is a quadratic residue
mod~$p$, $Q$ can also take the value
\[ Q\equiv\frac{p+1}{2} +\frac{p+1}{2} +2^2\equiv 5, \]
and repeating the argument sufficiently many times
(always with $R=2$) we either find an equation
for $k^2$ without any solution, or we can realize $Q\equiv p\equiv 0$ as
a sum of three squares mod~$p$. But then the equation
\[ Rk^2\equiv Rk^2-Q\equiv 0\]
does not have any solution $k$ coprime to $p$ if we
choose $R\not\equiv 0$ mod~$p$, as was
desired.  \hfill $\Box$
\section{Cyclic groups of prime power order}
We now show that ${\Bbb Z}_{p^k}$ is also a contact group, at least for
primes $p \geq 5$, by the same method as in the previous section.

\begin{prop}
\label{prop:power}
Write $h=h_{k,l}$ for the inclusion ${\Bbb Z}_{p^k}\rightarrow
{\Bbb Z}_{p^l}$, $k\leq l$. For $p\geq 5$ there is a short exact sequence
\[ 0\longrightarrow\om (B{\Bbb Z}_{p^{k-1}})
\stackrel{(Bh)_*}{\longrightarrow} \om (B{\Bbb Z}_{p^k})
\stackrel{(Bh)^!}{\longrightarrow} \om (B{\Bbb Z}_p)
\longrightarrow 0,\]
and $\om (B{\Bbb Z}_{p^k})$ is generated by lens spaces.
\end{prop}

\noindent {\em Proof.}
A spectral sequence argument as in the preceding section shows that
$\om (B{\Bbb Z}_{p^k})\cong \widetilde{\Omega}_5
(B{\Bbb Z}_{p^k})$ has order~$p^{2k}$.
The inclusion homomorphism $(Bh)_*$ is injective because the
corresponding homomorphism on homology is injective and the bordism
spectral sequence collapses at the $E^2$-page, cf.~\cite[(37.2)]{cofl64}.
Furthermore, the transfer homomorphism $(Bh)^!$ is surjective, for
we have shown that $\om (B{\Bbb Z}_p)$ is generated by ${\Bbb Z}_p$-lens
spaces, and every free linear ${\Bbb Z}_p$-action on $S^5$ extends to
a free linear ${\Bbb Z}_{p^k}$-action. Finally, the
composition 
\[ (Bh_{k-1,k})^!(Bh_{k-1,k})_*:\, \om (B{\Bbb Z}_{p^{k-1}})
\longrightarrow \om (B{\Bbb Z}_{p^{k-1}}) \]
is multiplication by~$p$, the index of ${\Bbb Z}_{p^{k-1}}$
in~${\Bbb Z}_{p^k}$
(see~\cite[(20.2)]{cofl64}, this is a general statement about the composition
of inclusion and transfer for central subgroups). Therefore the composition
\[ \om (B{\Bbb Z}_{p^{k-1}}) \longrightarrow \om (B{\Bbb Z}_{p^k})
\longrightarrow \om (B{\Bbb Z}_{p^{k-1}}) \longrightarrow
\om (B{\Bbb Z}_p) \]
is the zero map, because every
element in $\om (B{\Bbb Z}_p)$ has order~$p$ (here the
argument fails for~$p=3$). This proves that the sequence in the
proposition is exact, since the order of the middle group
is the product of the order of the two outer groups.

Arguing inductively, we assume that $\om (B{\Bbb Z}_{p^{k-1}})$ is
generated by lens spaces.
Given $u\in\om (B{\Bbb Z}_{p^k})$, we know that $(Bh)^!(u)\in
\om (B{\Bbb Z}_p)$ can be represented by a sum of ${\Bbb Z}_p$-lens spaces.
Lifting these ${\Bbb Z}_p$-actions to ${\Bbb Z}_{p^k}$-actions, we
get a sum $u_0$ of ${\Bbb Z}_{p^k}$-lens spaces such that
$(Bh)^!(u)=(Bh)^!(u_0)$. Notice, however, that the order of a
${\Bbb Z}_{p^k}$-lens space in $\om (B{\Bbb Z}_{p^k})$ is~$p^k$
\cite[(37.9)]{cofl64}, so the short exact sequence is not split.
Then $u-u_0=(Bh)_*(u_1)$ with $u_1$ represented by a sum of
lens spaces by the induction assumption.
This proves the proposition. \hfill $\Box$

\section{Proof of the reduction theorem}
As mentioned earlier, our proof of Theorem~\ref{thm:reduction}
differs from the corresponding proof in~\cite{kwsc90} only insofar as
we include some additional details, and that our situation is a bit
simpler because of the restriction to dimension five and
to odd order groups.

\vspace{2mm}

\noindent {\em Proof of Theorem}~\ref{thm:reduction}.
One direction of the theorem is the content of Lemma~\ref{lem:naturality}.
For the converse, we now assume that we are given $\alpha\in\Om$ with
$(Bj_p)^!\alpha\in\cont (\pi_p)$ for all primes $p$ dividing~$|\pi |$, and
we need to show that $\alpha\in\Cont$.

Write $T_p$ for the composition $(Bj_p)_*(Bj_p)^!$. While the composition
of inclusion and transfer (in this order) can be computed, at least
for normal subgroups (we used this in the proof of
Proposition~\ref{prop:power}), this is not true, in general, for a composition
of transfer and inclusion. We circumvent this problem by reducing the
computation of $T_p$ on bordism groups to that of the corresponding
homomorphism on homology groups.

First we reproduce an elementary algebraic lemma of~\cite{kwsc90}.

\begin{lem}
\label{lem:algebraic}
Let $R$ be a Noetherian ring, $\Omega$ a finitely generated $R$-module,
and $T$ an automorphism of $\Omega$. If $P$ is a submodule of $\Omega$ such
that $T(P)\subset P$, then $T(P)=P$.
\end{lem}

\noindent {\em Proof.}
The ascending chain of submodules
\[ P\subset T^{-1}(P)\subset T^{-2}(P)\subset ...\]
must terminate, since $R$ is Noetherian. Thus $T^{-m}(P)=T^{-m-1}(P)$
for some~$m$, which on applying $T^{m+1}$ yields $T(P)=P$.\hfill $\Box$

\vspace{2mm}

In the next lemma, ${\Bbb Z}_{(p)}$ denotes the integers localized
at~$p$ and $\Om_{(p)}$ the $p$-primary component of $\Om$.

\begin{lem}
\label{lem:p-primary}
For any prime $p$ dividing~$|\pi |$, the homomorphism $T_p\otimes
{\Bbb Z}_{(p)}$ is an isomorphism of $\Om_{(p)}$.
\end{lem}

\noindent {\em Proof.}
The cobordism spectral sequence yields the following commutative dia\-gram
with exact rows (except for the commutativity this follows from
the argument in the proof of Proposition~\ref{prop:cyclic},
since $H_*(B\pi )$ admits a $p$-primary decomposition with $p$
ranging over the primes dividing~$|\pi |$, cf.~\cite[III.10.2]{brow82}).

\[ \begin{array}{ccccc}
H_1(B\pi ) & \rightarrowtail & \Om & \twoheadrightarrow & H_5(B\pi )\\
\mbox{\Large $\downarrow$}  &  & \mbox{\Large $\downarrow$}  &
& \mbox{\Large $\downarrow$} \\
H_1(B\pi_p ) & \rightarrowtail & \om (B\pi_p) &
\twoheadrightarrow & H_5(B\pi_p )\\
\mbox{\Large $\downarrow$}  &  & \mbox{\Large $\downarrow$} &
& \mbox{\Large $\downarrow$}  \\
H_1(B\pi ) & \rightarrowtail & \Om & \twoheadrightarrow & H_5(B\pi )
\end{array} \]
The vertical arrows at the top denote the transfer homomorphism
$(Bj_p)^!$, those at the bottom the inclusion homomorphism $(Bj_p)_*$.
Commutativity of the squares on the right is proved in~\cite[(20.3)]{cofl64}.
Commutativity of the squares on the left follows similarly by
considering the isomorphism $\mu :\widetilde{\Omega}_1
(B\pi)\rightarrow H_1(B\pi )$ and the
inclusion of $\widetilde{\Omega}_1(B\pi )$
in $\widetilde{\Omega}_5(B\pi )$ by
tensoring with~$\Omega_4$ (and the same for~$\pi_p$). Alternatively,
this can be seen directly
from the geometric definitions of the maps in question.

On homology the composition $T_p=(Bj_p)_*(Bj_p)^!$ is multiplication
by the index of $\pi_p$ in~$\pi$ (cf.~\cite[III.9.5]{brow82}). Thus,
$T_p\otimes {\Bbb Z}_{(p)}$ is an isomorphism on homology localized at~$p$,
and by the five-lemma applied to the $p$-primary part of the diagram
above it is also an isomorphism on $\Om_{(p)}$.
This proves the lemma. \hfill $\Box$

\vspace{2mm}

By assumption we have $(Bj_p)^!\alpha\in\cont (\pi_p)$. Then by
Lemma~\ref{lem:naturality} we have $T_p\alpha\in\Cont$. So
$(T_p\otimes {\Bbb Z}_{(p)})(\alpha_{(p)})\in\Cont_{(p)}$,
and $\Cont_{(p)}$ is $(T_p\otimes {\Bbb Z}_{(p)})$-invariant by
that lemma. Then it follows from Lemmas~\ref{lem:algebraic}
and~\ref{lem:p-primary} that
$\alpha_{(p)}\in\Cont_{(p)}$.

Since $H_*(B\pi )$ admits a $p$-primary decomposition with $p$
ranging over the primes dividing~$|\pi |$,
the same holds for $\Om$, and so $\alpha_{(p)}\in\Cont_{(p)}$
for all $p$ dividing $|\pi |$ implies $\alpha\in\Cont$.
This concludes the proof of the reduction theorem.\hfill $\Box$

\vspace{2mm}

\noindent {\em Remark.}
Even though the bordism spectral sequence no longer collapses for
groups of even order, the reduction theorem still holds (by
essentially the same argument). Combining this with the fact
that ${\Bbb Z}_2$ is a contact group as proved in~\cite{geth98},
we see that in Theorem~1 we may actually allow that $|\pi |$ contains
a single prime factor~2.

\vspace{2mm}

In some instances one can be more specific about the contact manifolds which
generate $\Om$. To illustrate this, we briefly return to the
metacyclic groups of the introduction with $m=p^k$ for $p$ some prime
greater than or equal to five, $n=3$, and $r$ a primitive cube root
of 1 mod~$p^k$. Write $D_{p^k,3}$ for these groups.

As shown by Madsen~\cite[Theorem~4.13]{mads78}, there is a smooth
5-dimensional spherical space form $M_{p^k,3}$ with fundamental
group isomorphic to $D_{p^k,3}$ which is covered by lens spaces $L^5_3
\rightarrow M$ and $L^5_{p^k}\rightarrow M$ (Madsen's result is in fact
more general). In other words, with $\alpha\in\om (BD_{p^k,3})$
denoting the class of $M_{p^k,3}$ and the classifying map of its
universal covering, both $(Bj_3)^!\alpha$ and $(Bj_p)^!\alpha$ are
represented by lens spaces. By Theorem~\ref{thm:reduction}, $\alpha$
lies in $\cont (D_{p^k,3})$, and by Theorem~\ref{thm:cobord} we
know that $M_{p^k,3}$ admits a contact structure. By comparison,
Theorem~27 of~\cite{geig97} only guarantees the existence of a
contact structure on some special 5-dimensional space form with
fundamental group $D_{p^k,3}$, obtained via a construction of Petrie.

More can be said, however. The same spectral sequence argument as in
the proof of Proposition~\ref{prop:cyclic} shows that $\om (BD_{p^k,3})$
has order~$9p^k$. Now $(Bj_3)^!\alpha$ is represented by a ${\Bbb Z}_3$-lens
space and thus has order~9, whereas $(Bj_p)^!\alpha$ has order~$p^k$.
So $\alpha$ is an element of order (at least) $9p^k$
and therefore $\om (BD_{p^k,3})$ is a cyclic group ${\Bbb Z}_{9p^k}$,
generated by~$\alpha$. We have thus found a
5-dimensional spherical space form with fundamental group $D_{p^k,3}$
which carries a contact structure and generates~$\om (BD_{p^k,3})$.

\vspace{1mm}
\noindent
\begin{minipage}{5.5cm}
\noindent
{\em Mathematisch Instituut\\Universiteit Leiden\\Postbus
9512\\2300 RA Leiden\\The Netherlands\\e-mail: geiges@wi.leidenuniv.nl}
\end{minipage}
\begin{minipage}{6.5cm}
\noindent
{\em DPMMS\\University of Cambridge\\16 Mill Lane\\Cambridge
CB2 1SB\\United Kingdom\\e-mail: C.B.Thomas@dpmms.cam.ac.uk}
\end{minipage}
\end{document}